\documentclass[graybox]{svmult}
\usepackage{mathptmx}       
\usepackage{helvet}         
\usepackage{courier}        
\usepackage{type1cm}        
\usepackage{makeidx}         
\usepackage{graphicx}        
\usepackage{multicol}        
\usepackage[bottom]{footmisc}
\usepackage{subfigure}
\usepackage{url}

\usepackage{amsfonts}
\usepackage{amsmath}
\usepackage{amssymb}

\begin{document}

\title{Lippmann-Schwinger-Lanczos approach for inverse scattering problem of Schr\"{o}dinger equation in the resonance frequency domain}
\author{Anarzhan Abilgazy and Mikhail Zaslavskiy} 
\institute{Anarzhan Abilgazy \at Southern Methodist University, Department of Mathematics, Clements Hall, 3100 Dyer st,
Dallas, TX 75205  
\and Mikhail Zaslavskiy \at Southern Methodist University, Department of Mathematics, Clements Hall, 3100 Dyer st,
Dallas, TX 75205\at \email{mzaslavskiy@smu.edu} }  

\maketitle

\abstract{Reconstructions of potential in Schr\"{o}dinger equation with data in the diffusion frequency domain have been successfully obtained within Lippmann-Schwinger-Lanczos (LSL) approach in \cite{BaChDrMoZa}, however limited resolution away from the sensor positions resulted in rather blurry images. To improve the reconstructions, in this work we extended the applicability of the approach to the data in the resonance frequency domain. We proposed a specific data sampling according to Weyl's law that allows us to obtain sharp images without oversampling and overwhelming computational complexity. Numerical results presented at the end illustrate the performance of the algorithm.} 

\section{Introduction}
Lippmann-Schwinger (LS) approach has been shown to be a powerful tool for solving inverse scattering problems. It allows to recover medium properties using the frequency-domain or time-domain near-field measured data. The range of applications includes geophysical prospecting, radar and sonar imaging, medical imaging, deep space exploration, remote sensing, and many others. LS approach allows us to formulate these problems in terms of a nonlinear integral equation that involve both the unknown medium properties as well as the internal solution of the forward problem with those unknown properties. Born approximation linearizes that equation in the vicinity of the given background, however, it is known to be accurate only for small  perturbations of properties.

Reduced-order models (ROMs) have been successfully applied to solve inverse scattering problems in \cite{Borcea, druskin2016direct, druskin2013solution}. ROMs have been recently combined with Lippmann-Schwinger (LS) approach to improve the robustness of the latter \cite{DrMoZa, BaChDrMoZa}. This paper discusses an efficient data-driven ROM formulation we developed for the frequency-domain Schr\"{o}dinger equation to learn the internal solution via the measured data and background internal solution. When plugged into the LS equation, it enables a direct linear imaging of the medium properties without limitations of Born approximation. Such an approach has been developed for Schr\"{o}dinger equation in the diffusive frequency domain \cite{BaChDrMoZa}, however, due to severe ill-posedness of the inverse scattering problem \cite{Mandache:2001} it provides rather blurry images of the medium properties. Moreover, adding the measured data will not improve the image for the same reason. Here, we extended the applicability of the approach to the resonance frequency domain. In order to avoid overfitting and handling large datasets, we proposed an efficient data sampling for ROMs' interpolation based on the Weyl's law for asymptotic distribution of resonances.

\section{Schr\"{o}dinger equation}

We consider the 1D frequency-domain Schr\"{o}dinger equation with a potential $p(x)$ for a scalar function $u$ 
  \begin{equation} \label{eq:1}
 \begin{aligned}
    - \frac{d^2 u(x)}{dx^2} + p(x)u(x) + \lambda u(x) = \delta(x) 
    \quad \text{in} \quad \Omega = (0,L),
 \end{aligned}
\end{equation}
with Neumann boundary conditions $\left.\frac{du}{dx}\right|_{x=0} = \left.\frac{du}{dx}\right|_{x=L} =0$ . Here $0< L <  \infty$, and the source function $\delta(x)$ is a delta-function. The solution can be rewritten in terms of resolvent as $u(x,\lambda) = (\mathcal{L} + \lambda I)^{-1}\delta$ where $ \mathcal{L} = - \frac{d^2}{dx^2} + p I$. We define the single-input single-output (SISO) transfer function for $0>\lambda\notin sp(\mathcal{-L})\cup sp(\frac{d^2}{dx^2})$ as 
\begin{equation}\label{eq:2}
    F(\lambda) = \int_{0}^{L} \delta(x) u(x, \lambda) dx = \left< \delta, (\mathcal{L} + \lambda I)^{-1}\delta\right> =u(0,\lambda)
\end{equation}

\noindent where $\left<\cdot,\cdot\right>$ denote $L_2$-inner product. In the SISO inverse scattering problem we consider, the goal is to find $p(x)$ in (\ref{eq:1}) from the data that is given by

\begin{equation*}
\label{eq:data}
\begin{aligned}
        F(\lambda)|_{\lambda = \lambda_j} \in \mathbb{R}, \quad \frac{dF(\lambda)}{d \lambda}|_{\lambda = \lambda_j} \in \mathbb{R} \quad for \quad j =1,...,m.
\end{aligned}
\end{equation*}

\section{Lippmann-Schwinger approach}
 Consider background problem 
   \begin{equation} \label{eq:bckg}
 \begin{aligned}
    - \frac{d^2u^0(x)}{dx^2} + \lambda u^0(x) = \delta(x) 
    \quad \text{in} \quad \Omega,
 \end{aligned}
\end{equation}
with the boundary conditions     $\left.\frac{du^0}{dx}\right|_{x=0} =\left.\frac{du^0}{dx}\right|_{x=L}= 0.$
We can write similar equation as (\ref{eq:2}) for the transfer function $F_0$ of the background problem
 \begin{eqnarray}\label{eq:3}
     F_{0} (\lambda) = \int_{0}^{L} \delta(x) u^{0} (x, \lambda) dx=u^0(0,\lambda).
 \end{eqnarray}
Then using (\ref{eq:2}) and (\ref{eq:3}) we obtain Lippmann-Schwinger integral equation with respect to unknown coefficient $p(x)$
\begin{eqnarray}\label{eq:4}
    F_0(\lambda_j) - F(\lambda_j) = \int u^{0} (x,\lambda_j) u(x,\lambda_j) p(x) dx \quad for \quad j =1,...,m.
\end{eqnarray}
 Then after discretization of (\ref{eq:4}) using quadrature at nodes $\{x_i\}_{i = 1}^{n}$ we obtain the system of equations with respect to unknowns $\{p(x_i)\}^n_{i=1}$ 
\begin{equation*}
    F_0(\lambda_j) - F(\lambda_j) =  \left< u^0(x,\lambda_j), p u(x,\lambda_j) \right> \approx \sum^{n}_{i = 1}  h_i u^0 (x_i,\lambda_j) p(x_i)  u (x_i,\lambda_j),~j=1,\ldots, m
\end{equation*}
where $\{h_i\}^n_{i=1}$ are quadrature weights. This is a system of nonlinear equations with respect to $p(x)$ because $\{u(x,\lambda_j)\}^m_{j=1}$ themselves depend on $p(x)$. Born approximation replaces $u(x,\lambda_j)$ with $u^0(x,\lambda_j)$, however it is known to be accurate only for small $p(x)$. Below we will try to utilize the measured data $F$ in order to come up with a better approximant.

\section{Reduced Order Model}

Let $u_j$ be solutions to the equation (\ref{eq:1}) for $\lambda=\lambda_j,~j=1,...,m$. We consider ROM obtained via Galerkin formulation projecting the problem (\ref{eq:1}) onto the rational Krylov subspace ${\cal{V}} = colspan\{u_1 (x),...,u_m (x)\}$ \cite{Grimme}. Denoting $V=[u_1 (x),...,u_m (x)]\in \mathbb{R}^{\infty \times m}$, we obtain
 \begin{eqnarray} \label{eq:5}
     (S + \lambda M) c = b.
\end{eqnarray}
\noindent where in (\ref{eq:5}) $S$ is a stiffness matrix, $M$ is a mass matrix, and $b\in \mathbb{R}^m$ is column vector given by $b_i = \left< u_i, \delta \right> $. Both of stiffness and mass matrices are symmetric and positive definite matrices that are defined as $S_{ij} = \left< u_i , \mathcal{L} u_j  \right>,~ M_{ij} = \left< u_i ,  u_j  \right>$. We approximate the solution to the equation (\ref{eq:1}) by ROM as 
\begin{eqnarray}
\label{eq:galsol}
    u(x,\lambda) \approx \widehat{u}(x,\lambda) = V c = V \left( S + \lambda M \right)^{-1} b.
\end{eqnarray}

\noindent Thanks to the properties of Galerkin projection, the obtain that ROM satisfies the interpolation conditions $u(x,\lambda_j)=\widehat{u}(x,\lambda_j),~j=1,\ldots m$.  We also note that, though $\{u_i\}^n_{i=1}$ are not accessible, we can still compute the stiffness and mass matrices $S$ and $M$ in the data-driven way from using measured data (\ref{eq:data}) via the Loewner framework \cite{Antoulas} as

\begin{equation*}
    \begin{aligned}
        S_{ij} = \frac{\lambda_i F(\lambda_i) - \lambda_j F(\lambda_j)}{\lambda_i - \lambda_j}, \quad S_{ii} = \frac{ d(\lambda F)}{d \lambda} (\lambda_i);
    \end{aligned}
\end{equation*}

\noindent and 
\begin{equation*}
    \begin{aligned}
        M_{ij} = \frac{F(\lambda_j) - F(\lambda_i)}{\lambda_i - \lambda_j}, \quad M_{ii} = - \frac{dF}{d \lambda} (\lambda_i).
    \end{aligned}
\end{equation*}

As was first noticed in \cite{druskin2016direct} for the time-domain problem, the orthogonalized time-domain snapshots are weakly dependent on medium properties. Their frequency-domain counterparts can be obtained via Lanczos algorithm applied for the matrix pencil in (\ref{eq:5}) \cite{BaChDrMoZa}. This leads to tridiagonal matrix $T \in \mathbb{R}^{m \times m}$ and M-orthonormal Lanczos vectors $q_i \in \mathbb{R}^m$ such that for $Q=(q_1,\ldots,q_m)\in\mathbb{R}^{m\times m}$ we have
\begin{eqnarray*}
     SQ = QT, \quad Q^{T} M Q = I,
\end{eqnarray*}
Then Galerkin solution $\widehat{u}$ in (\ref{eq:galsol}) can be rewritten as 
\begin{eqnarray}
\label{eq:up}
    \widehat{u} (x,\lambda) = \sqrt{b^T M^{-1} b} V Q (T + \lambda I)^{-1} e_1. 
\end{eqnarray}

Following \cite{BaChDrMoZa} we use the approximate identity $V Q \approx V_0 Q_0$ and replace $V Q$ by $V_0 Q_0$ in (\ref{eq:up})
\begin{eqnarray*}
    \widehat{u} (x,\lambda)\approx \mathbf{u}(x,\lambda) = \sqrt{b^T M^{-1} b} V_0 Q_0 (T + \lambda I)^{-1} e_1.
\end{eqnarray*}
Finally, plugging the obtained approximants  $\mathbf{u}(x,\lambda_j),~j=1,\ldots, m$ in Lippmann-Schwinger equation instead of $u$ we end up with direct imaging algorithm via solving the system of linear equations
\begin{equation*}
    F_0(\lambda_j) - F(\lambda_j) = \sum^{n}_{i = 1}  h_i u^0 (x_i,\lambda_j) p(x_i)  \mathbf{u} (x_i,\lambda_j),~j=1,\ldots, m
\end{equation*}

We note that the choice of data sampling points $\{\lambda_j\}^m_{j=1}$ is crucial for the efficiency of the approach. Indeed, undersampled data may ruin the image quality while oversampling results in a severe ill-posedness and an increase of the computational complexity. We found that sampling the data in accordance with Weyl's law provides accurate images at a reasonable computational cost. In particular, Weyl's law asymptotically counts the number $N(\lambda)$ of Neumann eigenvalues (counting multiplicities) that are less than or equal to $\lambda$ as $\lim_{\lambda\to\infty} \frac{N(\lambda)}{\lambda^d} = (2 \pi)^{-d} w_d \textbf{Vol}(\Omega)$ where $w_d$ is volume of unit ball in $\mathbb{R}^d$. In our 1D numerical experiments we chose three to five sampling points between every two approximate resonances given by Weyl's law. We note that such choice of $\{\lambda_j\}^m_{j=1}$ may still results in overfitting and, consequently, may require a regularization \cite{BaChDrMoZa}, however ill-posedness will be rather mild.

\section{Numerical Results}

\begin{figure}
    \centering
    \subfigure[]{\includegraphics[width=0.24\textwidth]{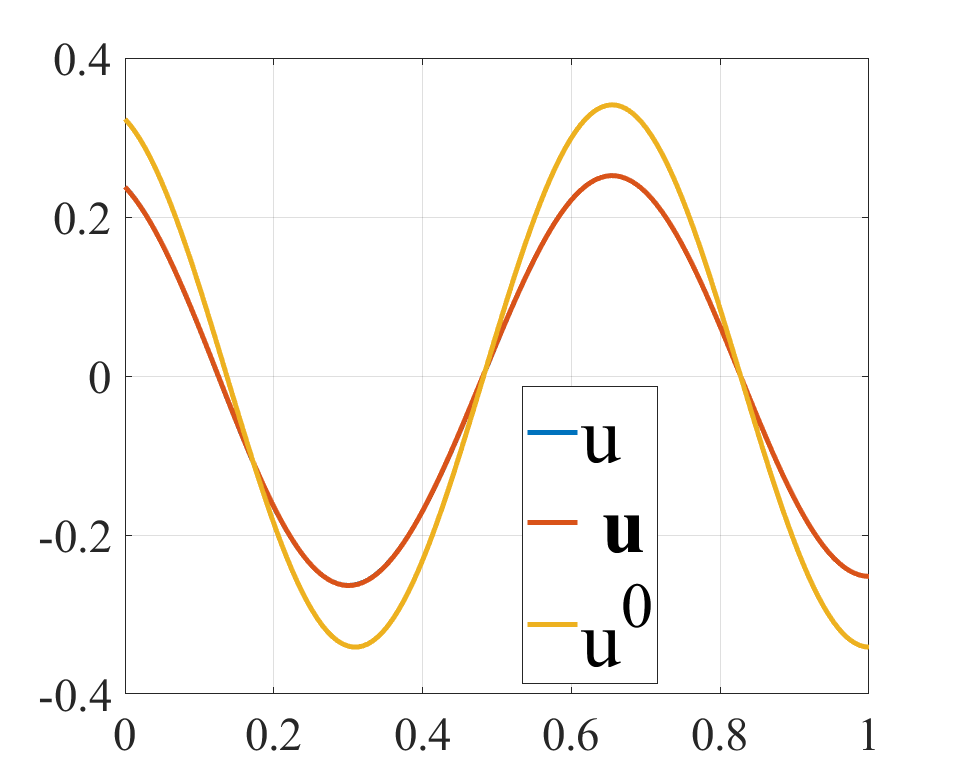}} 
    \subfigure[]{\includegraphics[width=0.24\textwidth]{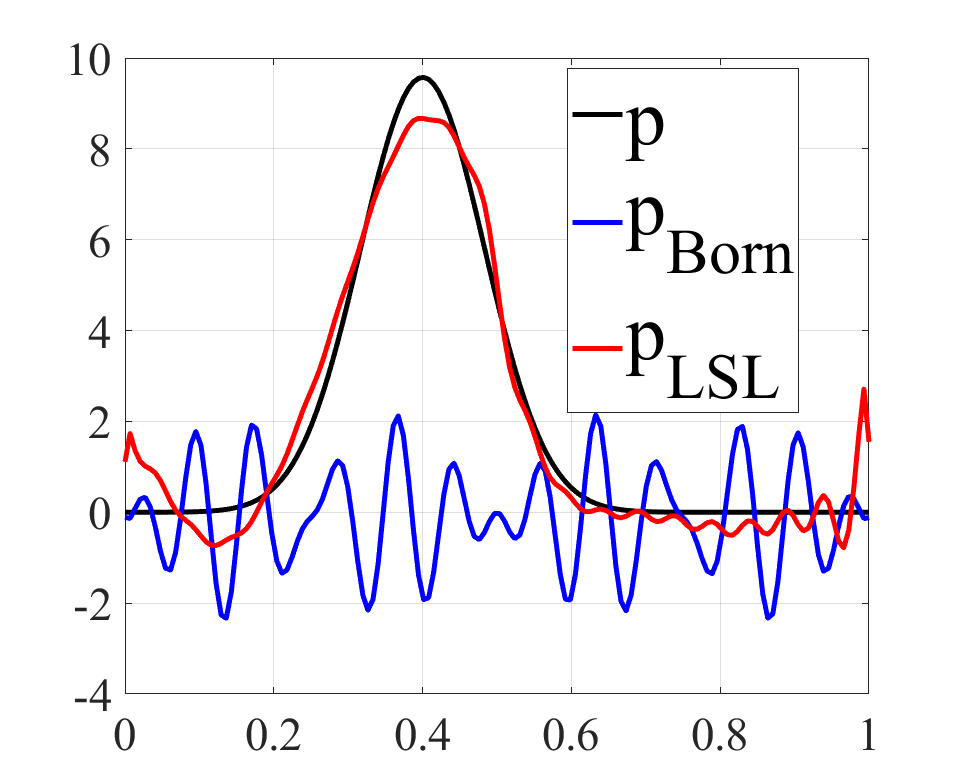}} 
    \subfigure[]{\includegraphics[width=0.24\textwidth]{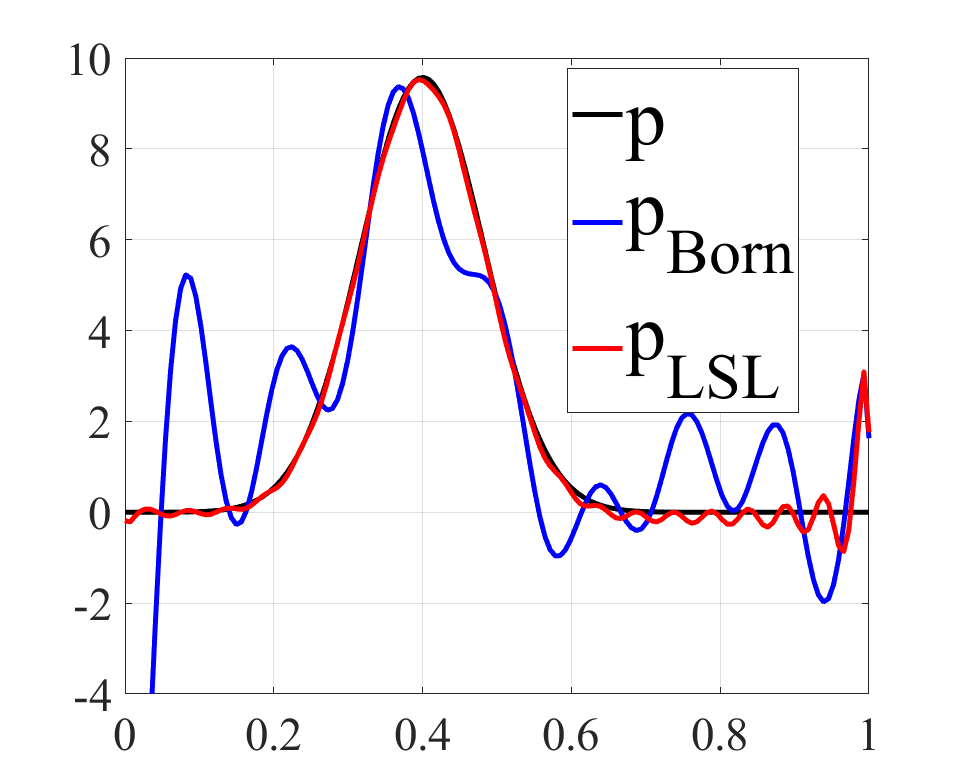}} 
    \subfigure[]{\includegraphics[width=0.24\textwidth]{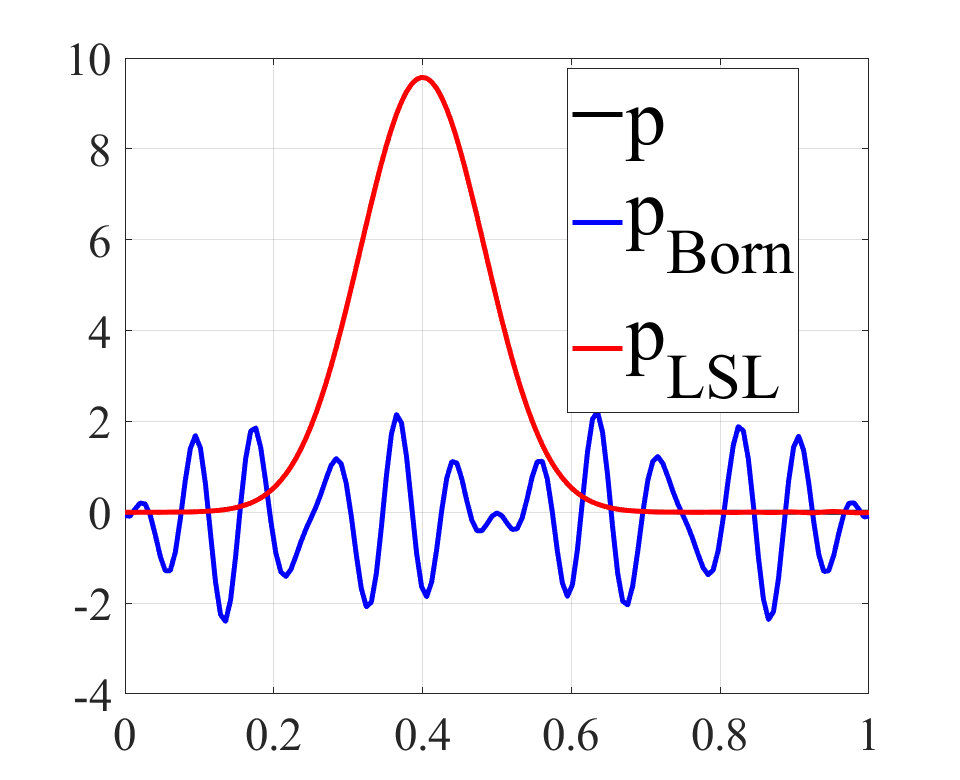}}
    \caption{$\mathbf{u}$ and $u$ almost coincide while $u^0$ is totally different (see (a) plot for $f=4$). True continuous potential $p$, its reconstructions $p_{Born}$ via Born and $p_{LSL}$ via LSL for $f=3$ see (b), $f=4$ (see (c) and $f=5$ (see (d))}
    \label{fig:ps}
\end{figure}

\noindent In our 1D experiments with noiseless data we took $N(\lambda) = 10$ and we compared results for $f=3,~4,~5$ equally spaced $\lambda$'s between every two approximate resonances. In the first experiment we considered the reconstruction of smooth potential $p(x)$ of Gaussian shape (see black curves in Fig. \ref{fig:ps} (b)-(d)). In Fig. \ref{fig:ps} (a) we plotted true solution $u(x,\lambda)$, its background counterpart $u^0(x,\lambda)$ as well as ROM solution $\mathbf{u}(x,\lambda)$ for the case $f=4$ for some intermediate value of $\lambda$. As one can observe, though the background solution $u^0$ looks totally different from $u$, the ROM approach still managed to reconstruct the internal solution $u$ rather accurately. That, in turn, results in improved imaging results (see Fig. (\ref{fig:ps}) (b), (c), (d) for the cases $f=3,~4,~5$, respectively). In particular, for $f=5$ the reconstructed potential almost coincides with the true one.

\begin{figure}
    \centering
    \subfigure[]{\includegraphics[width=0.3\textwidth]{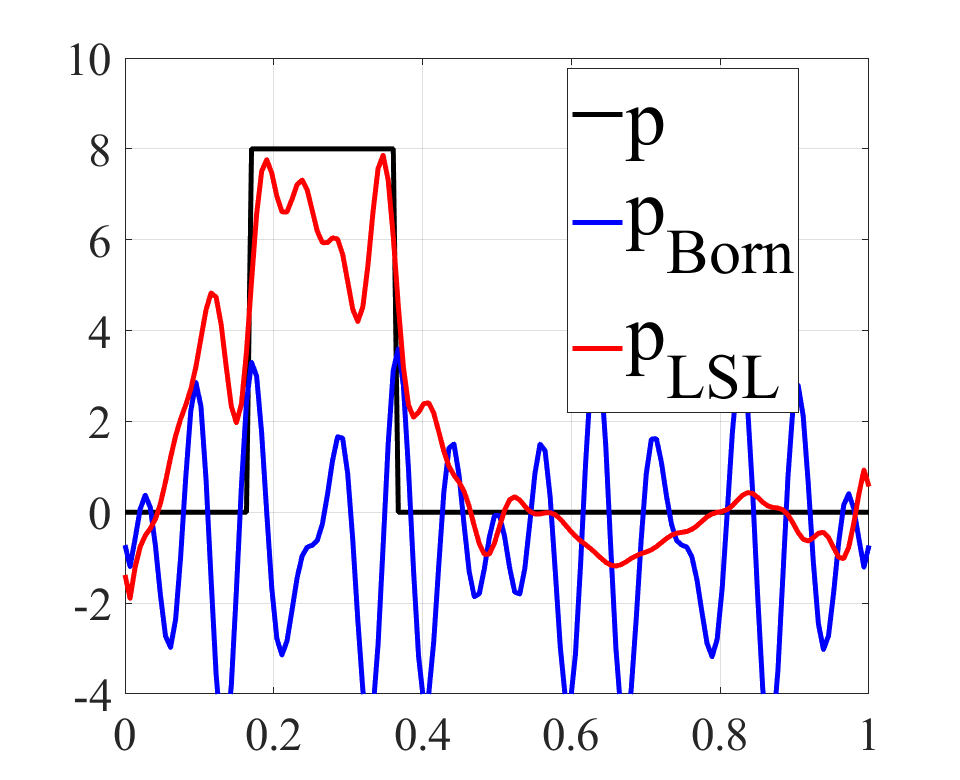}} 
    \subfigure[]{\includegraphics[width=0.3\textwidth]{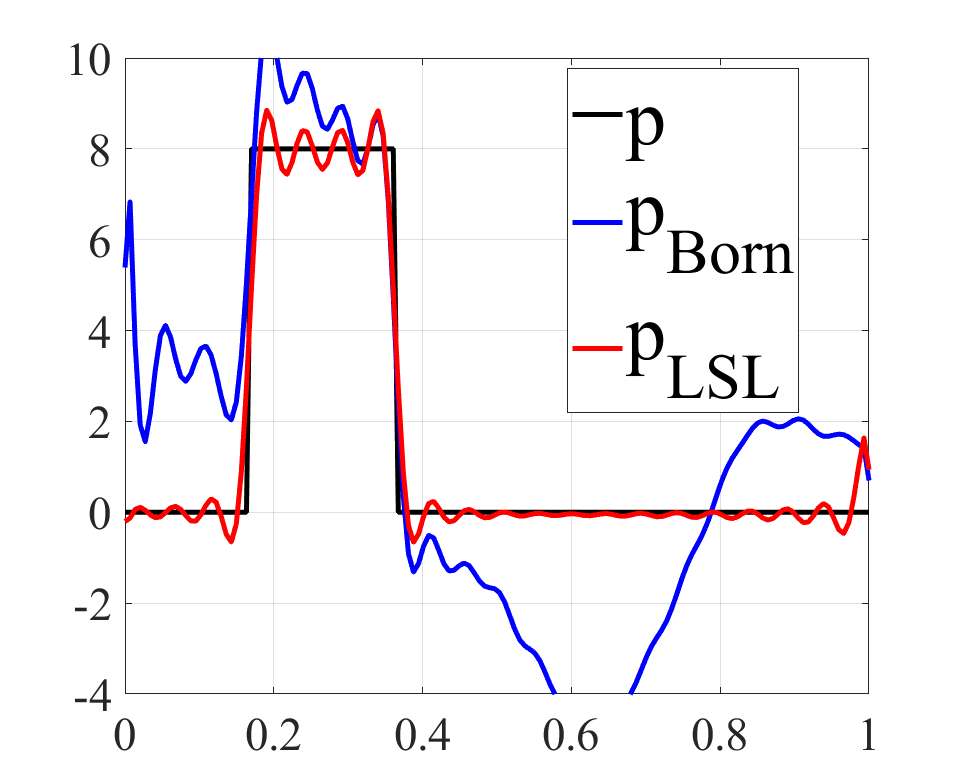}} 
    \subfigure[]{\includegraphics[width=0.3\textwidth]{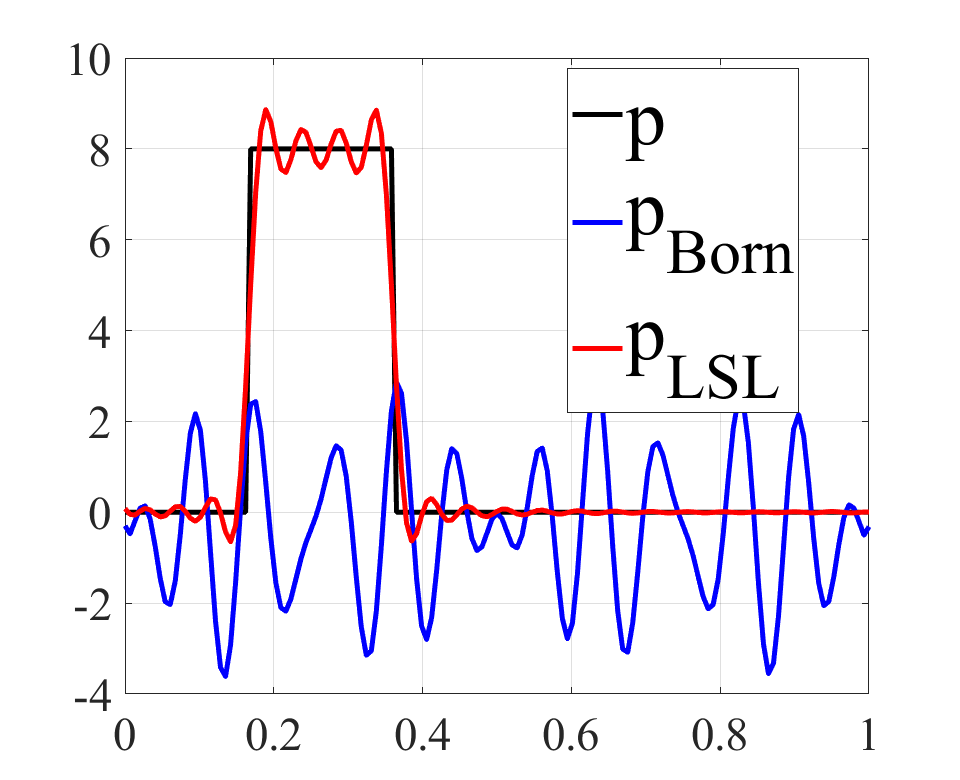}}
    \caption{True discontinuous potential $p$, its reconstructions $p_{Born}$ via Born and $p_{LSL}$ via LSL} 
    \label{fig:ps_step}
\end{figure}

\noindent In our second experiment we considered discontinuous potential $p(x)$ (see black curve in Fig. \ref{fig:ps_step}). Similar to the previous example, LSL images are more accurate compared to Born ones for all the cases $f=3,~4,~5$, though the results are polluted by Gibbs effects. 

\section{Conclusion}
In this work we extended the applicability of the LSL approach to the inverse scattering problem for the Schr\"{o}dinger equation in resonance frequency domain. That allowed to obtain a direct imaging algorithm that produced sharpened reconstructions compared to the diffusive frequency domain. To avoid the increase of computational cost and overfitting we sampled the data in accordance with Weyl's law putting several points between every two asymptotic resonances. We note that in higher dimensions ($d>1$) where the inverse scattering problem is over-determined it may be enough to use two sampling points between resonances, however in our 1D case we had to use three and more points. In our future work we plan to investigate multi-dimensional scenarios deeper and to consider noisy datasets.

\begin{acknowledgement}
The authors are  grateful to Justin Baker, Elena Cherkaev, Vladimir Druskin and Shari Moskow for productive discussions that inspired this research.
M. Zaslavskiy was partially supported by AFOSR grant  FA9550-20-1-0079.
\end{acknowledgement}




\end{document}